\documentclass[a4paper,11pt,english]{article}

\usepackage{amssymb}
\usepackage{mathtools}    
\usepackage[all]{xy}
\usepackage{amsthm}
\usepackage{currvita}
\usepackage{mathdots}
\usepackage{amsmath}
\usepackage{verbatim}
\usepackage{MnSymbol}
\usepackage{makeidx}
\usepackage{mathrsfs}
\usepackage[ansinew]{inputenc}
\usepackage{microtype}
\SelectTips{cm}{12}
\makeindex
\newtheorem{thm}{Theorem}[section]
\newtheorem{lemma}[thm]{Lemma}
\newtheorem{defi}[thm]{Definition}
\newtheorem{cor}[thm]{Corollary}

\newtheorem{prop}[thm]{Proposition}

\DeclarePairedDelimiter{\norm}{\lVert}{\rVert}
\renewcommand{\ker}{\operatorname{ker}}

\newcommand{\id}{\operatorname{id}}

\newcommand{\lin}{\operatorname{lin}}
\newcommand{\rank}{\operatorname{rank}}

\newcommand{\sgn}{\operatorname{sgn}}
\newcommand{\vp}{{\varphi}}

\newcommand{\Epi}{\operatorname{Epi}}

\newcommand{\TRO}{\operatorname{TRO}}

\newcommand{\Tri}{\operatorname{Tri}}

\newcommand{\ra}{{\mathbf R}}

\newcommand{\C}{{\mathbb C}}
\newcommand{\N}{{\mathbb N}}

\newcommand{\lr}[1] {{\left(#1\right)}}
\newcommand{\fr}[1] {{\frac{1}{#1}}}

\newcommand{\jb}{\text{$JB^*$-triple system}}

\newcommand{\bt}{\begin{thm}}
\newcommand{\et}{\end{thm}}
\newcommand{\bp}{\begin{proof}}
\newcommand{\ep}{\end{proof}}
\newcommand{\bl}{\begin{lemma}}
\newcommand{\el}{\end{lemma}}
\newcommand{\bc}{\begin{cor}}
\newcommand{\ec}{\end{cor}}
\newcommand{\bd}{\begin{description}}
\newcommand{\ed}{\end{description}}

\newcommand{\fA}{{\mathfrak A}}

\newcommand{\fG}{{\mathfrak G}}

\newcommand{\M}{\mathbb M}

\newcommand{\sub}{\subseteq}

\begin{document}
\thispagestyle{empty}
\begin{center}
\huge
\textsf{The universal enveloping TRO of a $JB^*$-triple system}\\
\vspace*{1cm} \large
Dennis Bohle, Wend Werner\\
Fachbereich Mathematik und Informatik\\
Westf\"alische Wilhelms-Universit\"at\\
Einsteinstra\ss e 62\\
48149 M\"unster \\
\vspace*{0.5cm}
 dennis.bohle@math.uni-muenster.de,\\ wwerner@math.uni-muenster.de\\
\vspace*{0.3cm}
\end{center}
\begin{abstract}
We associate to every $JB^*$-triple system a so called universal
enveloping TRO. We compute the universal enveloping TROs of the
finite dimensional Cartan factors.

 \footnotetext{The first
author was supported by the Graduiertenkolleg f\"ur analytische
Topologie und Metageometrie}

\end{abstract}

The content of the following is part of a project of the authors that has
as its objective to show that Cartan's classification of the (hermitian)
symmetric
spaces has a $K$-theoretic background. This project will come to a
conclusion in
the follow-up \cite{BoWe2}.

The symmetric spaces that will come up here consist of the open unit
balls of
so called $JB^*$-triples, an important generalization of the concept of a
$C^*$-algebra.
If the dimension is finite, their open unit balls coincide exactly with
the hermitian
symmetric spaces of non-compact type so that all of these spaces are
obtained trough
duality.

In the present paper we overcome a difficulty that is one of the main
obstacles for a direct generalization
of the K-theory of $C^*$-algebras: The impossibility, in general, to define
tensor products of a $JB^*$-triple
with the $n\times n$ matrices over the complex numbers. The $JB^*$-triples
that do have this property
are precisely the ternary rings of operators which as spaces coincide
with the
class of (full) Hilbert-$C^*$-modules.

We will study a construction which allows the passage from an arbitrary
$JB^*$-triple to such
a ternary ring in a way, that behaves so nicely that it will pave the way
for the program ahead.
In the first section we will collect some definitions and preliminary
results, the second
sections contains the actual construction of the enveloping ternary ring
of operators, whereas in
the third section the enveloping ternary rings of all finite-dimensional
Cartan factors are calculated. We do this quite differently from the
approach in \cite{BunceFeelyTimoneyI}
(the results of which were roughly obtained around the same time) in
that we use \emph{grids}. These objects will be helpful in the sequel
paper \cite{BoWe2} and are reminiscent of the root systems that are
central to the classical approach. Finally, in section 4, we slightly
improve a result from \cite{BunceFeelyTimoneyI} on the structure of the
enveloping ternary ring in some special cases.

All the results in the present paper are taken from the first named
author's PhD-thesis
\cite{Bohle}.

\section{Preliminaries}
We will first provide some notation, definitions and well known facts of triple theory. Our general references for the theory of $JB^*$-triple systems are \cite{Isidro-AglimpseatthetheoryofJordanBanachtriplesystems} and \cite{Upmeier-SymmetricBanachmanifoldandJordaCalgebras}. For $n\in\N$ we denote by $\M_n$ the $n\times n$-matrices over the complex numbers and if $Z$ is a Banach space, then $B(Z)$ is the Banach algebra of bounded linear operators on $Z$.
A Banach space $Z$ together with a sesquilinear mapping
\[Z\times Z\ni (x,y)\mapsto x\Box y\in B(Z) \]
is called a $\mathbf{JB^*}$-\textbf{triple
system}\index{$JB^*$-triple system}, if for the triple product
\[\{x,y,z\}:=(x\Box y)(z)\]
and all $a,b,x,y,z\in Z$ the following conditions are fulfilled:

The triple product $\{x,y,z\}$ is continuous in $(x,y,z)$, it is symmetric in the outer variables and
the $C^*$-condition $\norm{\{x,x,x\}}=\norm{x}^3$ is fulfilled.
The \textbf{Jordan triple identity}\index{Jordan triple identity} $$ \{a,b,\{x,y,z\}\}=\{\{a,b,x\},y,z\}-\{x,\{b,a,y\},z\}+\{x,y,\{a,b,z\}\}$$ holds.
The operator $x\Box x$ has non-negative spectrum in the Banach algebra $B(Z)$ and it is hermitian (i.e.\ $\exp(it(x\Box x))$ is isometric for all $t\in\mathbb{R}$).

A closed subspace $W$ of a $\jb$ $Z$ which is invariant under the
triple product, and therefore is a $JB^*$-triple system itself, is
called a $\mathbf{JB^*}$\textbf{-subtriple}\index{$JB^*$-subtriple} (or subtriple for short)
of $Z$.

 A closed subspace $I$ of a $\jb$ $Z$ is called a
$\mathbf{JB^*}$\textbf{-triple} \textbf{ideal}\index{ideal of a
$JB^*$-triple system}, if
$
\{Z,I,Z\}+\{I,Z,Z\}\subseteq I.
$
$JB^*$-triple ideals of $Z$ are $JB^*$-subtriples of $Z$ and the
kernel of a $JB^*$-triple homomorphism is always a $JB^*$-triple ideal.

Every $C^*$-algebra $\fA$
becomes a $JB^*$-triple system under the product
\begin{equation}\label{symmetrisiertes Product}
\{a,b,c\}:=\fr2\lr{ab^*c+cb^*a}.
\end{equation}
More generally every closed subspace of a $C^*$-algebra which is invariant under the product (\ref{symmetrisiertes Product}) is a $JB^*$-triple, called $\mathbf{JC^*}$-\textbf{triple system}.
A $JB^*$-triple system $Z$ which is a dual Banach space is called a
$\mathbf{JBW^*}$-\textbf{triple system}. Its predual is usually denoted by $Z_*$. The triple product
of a $JBW^*$-triple is separately $\sigma(Z,Z_*)$-continuous and its
predual
is unique.

An important example of $JB^*$-triples is given by the \textbf{TRO}s (ternary rings of operators). These are closed subspaces $T\sub B(H)$ such that \begin{equation}\label{TROprod}xy^*z\in T\end{equation} for all $x,y,z\in T$. \textbf{SubTRO}s are closed subspaces $U\sub T$ closed under (\ref{TROprod}) and \textbf{TRO-ideal}s are subTROs $I$ of $T$ such that $IT^*T+TI^*T+TT^*I\sub I$. TROs become $JB^*$-triples under the product (\ref{symmetrisiertes Product}).

Let $Z$ be a $JB^*$-triple system. An element $e\in Z$ that satisfies $\{e,e,e\}=e$ is called a \textbf{tripotent}. The collection of all non-zero tripotents in $Z$ is denoted by $\Tri(Z)$. A tripotent is called \textbf{minimal} if $\{e,Z,e\}=\mathbb{C}e$.
If $e$ is a non-zero tripotent, then $e$ induces a decomposition of $Z$ into the eigenspaces of $e\Box e$, the \textbf{Peirce decomposition}\[Z=P^e_0(Z)\oplus P^e_1(Z)\oplus P^e_2(Z),\] where $P^e_k(Z):=\{z\in Z:\{e,e,z\}=\frac{k}{2}z\}$ is the $\frac{k}{2}$-eigenspace, the \textbf{Peirce-k-space}, of $e\Box e$, for $k=0,1,2$.
Each Peirce-$k$-space, $k=0,1,2$, is again a $JB^*$-triple system.
In the case of a TRO $T$, the Peirce-2-space $P^e_2(T)$ becomes a unital $C^*$-algebra under the product $a\bullet b:=ae^*b$, denoted by $P^e_2(T)^{(e)}$.

Every finite-dimensional $JB^*$-triple system is the direct sum of so called Cartan factors $\mathcal{C}_1,\ldots,\mathcal{C}_6$. The two exceptional Cartan factors $\mathcal{C}_5$ and $\mathcal{C}_6$ can be realized as subspaces of the $3\times3$-matrices over the complex Cayley algebra $\mathbb{O}$ and the other four types are treated in detail in Section $4$.

\section{Universal objects}\label{section universal enveloping TRO}

We prove the existence of the universal enveloping TRO and the universal enveloping $C^*$-algebra of a $JB^*$-triple system. As a corollary we obtain a new proof of one of the main theorems of $JB^*$-triple theory.

The following lemma and theorem are generalizations of classical results for real
$JB$-algebras (cf.\ \cite{HancheOlsenStoermerJordanoperatoralgebras},
Theorem 7.1.3 and
\cite{alfsenschultz-Geometryofstatespacesofoperatoralgebras}, Theorem
4.36).

\begin{lemma}\label{universellerHilbertraum}
Let $Z$ be a $JB^*$-triple system. Then there exists a Hilbert space
$H$ such that for every $JB^*$-triple homomorphism $\varphi: Z\to
B(K) $ the $C^*$-algebra $\mathfrak{A}_{\varphi}$ generated by
$\varphi(Z)$ can be embedded $*$-isomorphically into $B(H)$.
\begin{proof}
The cardinality of $\varphi(Z)$ is less or equal to the cardinality
of $Z$. One can now proceed similar to the proof of \cite{alfsenschultz-Geometryofstatespacesofoperatoralgebras}, Lemma
4.35.
\end{proof}
\end{lemma}

\begin{thm}\label{universalTRO}
Let $Z$ be a $JB^*$-triple system.
\begin{description}
\item[(a)] There exist up to $*$-isomorphism a unique $C^*$-algebra
$C^*(Z)$ and a $JB^*$-triple homomorphism $\psi_Z:Z\to C^*(Z)$ such
that
\begin{description}
\item[(i)] For every $JB^*$-triple homomorphism $\varphi:Z\to
\mathfrak{A}$, where $\mathfrak{A}$ is an arbitrary $C^*$-algebra, exists a
$*$-homomorphism $C^*(\varphi): C^*(Z)\to \mathfrak{A}$ with
$C^*(\varphi)\circ\psi_Z=\varphi$.
\item[(ii)] $C^*(Z)$ is generated as a $C^*$-algebra by $\psi_Z(Z)$.

\end{description}
\item[(b)] There exists up to TRO-isomorphism a unique TRO
$T^*(Z)$ and a $JB^*$-triple homomorphism $\rho_Z:Z\to T^*(Z)$ such
that
\begin{description}
\item[(i)] For every $JB^*$-triple homomorphism $\alpha:Z\to
T$, where $T$ is an arbitrary TRO, exists a TRO-homomorphism $T^*(\alpha):
T^*(Z)\to T$ with $T^*(\alpha)\circ\rho_Z=\alpha$.
\item[(ii)] $T^*(Z)$ is generated as a TRO by $\rho_Z(Z)$.
\end{description}
\end{description}
\begin{proof}
Let $H$ be the Hilbert space from Lemma
\ref{universellerHilbertraum} and $I$ the family of $JB^*$-triple
homomorphisms from $Z$ to $B(H)$. Let
$\psi_Z:=\rho_Z:=\bigoplus_{\psi\in I}\psi$ and
$\hat{H}:=\bigoplus_{\psi\in I}H_\psi$ be $l^2$-direct sums with
$H_\psi:=H$. Then $\psi_Z$ and $\rho_Z$ are $JB^*$-triple
homomorphisms from $Z$ to $B(\hat{H})$. Let $C^*(Z)$ be the
$C^*$-algebra and $T^*(Z)$ the TRO generated by $\rho(Z)$ in
$B(\hat{H})$. If $\mathfrak{A}$ is a $C^*$-algebra and $\varphi:Z\to
\mathfrak{A}$ is a $JB^*$-triple homomorphism, where $\varphi(Z)$
w.l.o.g.\ generates $\mathfrak{A}$ as a $C^*$-algebra, then we can
suppose that $\mathfrak{A}$ is a subset of $B(H)$. Therefore
$\varphi$ can be regarded as an element of $I$. Let
$\pi_\varphi:\bigoplus_{\psi\in I}B(H_\psi)\to B(H_\varphi)$ be the
projection onto the $\varphi$-component, then
$\pi_\varphi(\psi_Z(z))=\pi_\varphi(\rho_Z(z))=\varphi(z)$ for all
$z\in Z$. We define $C^*(\varphi)$ resp. $T^*(\varphi)$ to be the
restrictions of $\pi_\varphi$ to $C^*(Z)$ resp. $T^*(Z)$. Uniqueness
is proved in the usual way using the universal properties.
\end{proof}
\end{thm}

We call $(T^*(Z),\rho_Z)$ the \textbf{universal enveloping TRO}\index{universal enveloping TRO} and
$(C^*(Z),\psi_Z)$  the \textbf{universal enveloping} $\mathbf{C^*}$\textbf{-algebra}\index{universal enveloping $C^*$-algebra} of $Z$
respectively. Most of the time we only use $T^*(Z)$ and $C^*(Z)$
as shorter versions.

Similar to the classical case \cite{alfsenschultz-Geometryofstatespacesofoperatoralgebras}, Proposition 4.40 there exists a TRO-antiautomorphism on $T^*(Z)$:

\begin{prop}\label{canonical TRO-antiautomorphism}
Let $Z$ be a $JB^*$-triple system.
There exists a TRO-antiautomorphism $\theta$ (i.e.\ a linear, bijective mapping from $T^*(Z)$ to $T^*(Z)$ such $\theta(xy^*z)=\theta(z)\theta(y)^*\theta(x)$ for all $x,y,z\in T^*(Z)$) of $T^*(Z)$  of order $2$ such that $\theta\circ\rho_Z=\rho_Z.$
\begin{proof}
Denote by $T^*(Z)^{op}$ the opposite TRO of $T^*(Z)$, i.e.\ the TRO that coincides with $T^*(Z)$ as a set and is equipped with the same norm, if $\gamma: T^*(Z)\to T^*(Z)^{op}$, $\gamma(a)=a^{op}$ is the identity mapping then $\left(xy^*z\right)^{op}=z^{op}(y^{op})^*x^{op}$ for all $x,y,z\in T^*(Z)$.

The composed mapping $\gamma\circ\rho_Z:Z\to T^*(Z)^{op}$ is a $JB^*$-triple homomorphism and thus lifts to a TRO-homomorphism $T^*(\gamma\circ\rho_Z):T^*(Z)\to T^*(Z)^{op}$. We put $$\theta:=\gamma^{-1}\circ T^*(\gamma\circ\rho_Z):T^*(Z)\to T^*(Z).$$ Easily it can be seen (since $\theta$ fixes by construction $\rho_Z(Z)$ which generates $T^*(Z)$ as a TRO) using the universal properties of $T^*(Z)$ that $\theta$ is a TRO-anti\-auto\-morphism of order $2$.
\end{proof}
\end{prop}

We refer to $\theta$ as the \textbf{canonical TRO-antiautomorphism of order} $\mathbf{2}$\index{canonical TRO-antiautomorphism of order $2$} on $T^*(Z)$.

\begin{cor}
If the $JB^*$-triple system $Z$ in Theorem \ref{universalTRO} is a
$JC^*$-triple then the mappings $\psi_Z$ and $\rho_Z$ are injective.
\end{cor}

Obviously $\psi_Z$ and $\rho_Z$ are the $0$ mappings if $Z$ is
purely exceptional.

\begin{lemma}\label{FortsetzenVonHom}
For every $JB^*$-triple ideal $I$ in a $JB^*$-triple system $Z$ and
every $JB^*$-triple homomorphism $\varphi:I\to W$, where $W$ is a
$JBW^*$-triple system, there exists a $JB^*$-triple homomorphism
$\Phi:Z\to W$ which extends $\varphi$.
\begin{proof}
We know by
\cite{Dineen-CompleteholomorphicvectorfieldsontheseconddualofaBanachspace}
that the second dual $Z''$ of $Z$ is a $JBW^*$-triple system and the
canonical embedding $\iota:Z\to Z''$ is an isometric $JB^*$-triple
isomorphism onto a norm closed $w^*$-dense subtriple of $Z''$. By
\cite{BunceandChu-CompactoperationsmultipliersandRadonNikodympropertyinJBtriples}, Remark 1.1
and since $W$ is a $JBW^*$-triple system there exists a unique,
$w^*$-continuous extension $\overline{\varphi}:I''\to W$ of
$\varphi$ with
$\overline{\varphi}(I'')=\overline{\varphi(I)}^{w^*}$. Let
\[I^\perp:=\{x\in Z'':y\mapsto \{x,i,y\}\text{ is the }0\text{
mapping for all }i\in I''\}\] be the $w^*$-closed orthogonal
complement of $I''$ with $Z''=I''\oplus I^\perp$ (cf.\
\cite{Horn-CharacterizationofthepredualandidealstructureofaJBW-triple},
Theorem 4.2 (4)). If we denote the projection of $Z''$ onto $I''$
with $\pi$ we get the desired extension of $\varphi$ by defining $\Phi:=\overline{\varphi}\circ\pi\circ\iota$.
\end{proof}
\end{lemma}

We obtain a new proof of an important theorem of Friedman and Russo
(cf.\ \cite{FriedmanRusso-TheGelfandNaimarktheoremforJBtriples},
Theorem 2):
\begin{cor}
Any $JB^*$-triple system $Z$ contains a unique purely exceptional
ideal $J$ such that $Z/J$ is $JB^*$-triple isomorphic to a
$JC^*$-triple system.
\begin{proof}
Let $J$ be the kernel of the mapping $\rho_Z:Z\to T^*(Z)$, which is
a $JB^*$-triple ideal. We know that $Z/J$ is a
$JB^*$-triple system which is $JB^*$-triple isomorphic to the
$JB^*$-triple system $\rho_Z(Z)\subseteq T^*(Z)$ and hence to a
$JC^*$-triple system.

Let us assume that $J$ is not purely exceptional which means that
there exists a non-zero $JB^*$-triple homomorphism $\varphi$ from $J$
into some $B(H)$. This $JB^*$-triple homomorphism extends by Lemma \ref{FortsetzenVonHom} to a $JB^*$-triple homomorphism $\phi:Z\to B(H)$.
Since $\phi=T^*(\phi)\circ \rho_Z$ holds, $\phi$ vanishes on $J$,
which is a contradiction.

Now let $I$ be another purely exceptional ideal such that $Z/I$ is
$JB^*$-triple isomorphic to a $JC^*$-triple system. On the one hand we
have $I\subseteq\ker(\rho_Z)=J$. On the other hand let $\varphi:Z\to
B(H)$ be a $JB^*$-triple homomorphism with kernel $I$. Then
$\varphi$ has to vanish on $J$ and therefore $J\subseteq I$.
\end{proof}
\end{cor}

\section{Cartan factors}\label{section Cartan factors}

In this section we compute
the universal enveloping TROs of the finite-dimensional Cartan
factors. Since the universal enveloping TROs of the two exceptional
factors are $0$, we have to compute the factors of type I--IV. We
do so by using the grids spanning these factors (cf.\
\cite{dangfriedman-Classificationoftriplefactorsandapplications} and Chapter 2). We make heavy use of the elaborate work on grids in \cite{nealrusso-Contractiveprojectionsandoperatorspaces}.

\subsection{Factors of type IV}\label{Factors of type IV}

A \textbf{spin system} is a subset $S=\{\id,s_1,\ldots,s_n\}$, $n\geq 2$, of
self-adjoint elements of $B(H)$ which satisfy the anti-commutator
relation $s_is_j+s_js_i=2\delta_{i,j}$ for all $i,j\in\{1,\ldots,n\}$. The complex linear span of
$S$ is a $JC^*$-algebra of dimension $n+1$ (cf.\
\cite{HancheOlsenStoermerJordanoperatoralgebras}). Every
$JC^*$-triple system which is $JB^*$-isomorphic to such a $JC^*$-algebra is
called a \textbf{spin factor}. We now recall the definition of a spin grid: A
\textbf{spin grid} is a collection $\{u_j,\widetilde{u}_j|j\in J\}$ (or
$\{u_j,\widetilde{u}_j|j\in J\}\cup \{u_0\}$ in finite odd
dimensions), where $J$ is an index set with $0\notin J$, for $j\in
J$, $u_j,\widetilde{u}_j$ are minimal tripotents and, if we let
$i,j\in J$, $i\neq j$, then
\begin{description}
\item[(SPG1)] $\{u_i,u_i,\widetilde{u}_j\}=\frac{1}{2}\widetilde{u}_j$,
$\{\widetilde{u}_j,\widetilde{u}_j,u_i\}=\frac{1}{2}u_i$,
\item[(SPG2)] $\{u_i,u_i,u_j\}=\frac{1}{2}u_j$,
$\{u_j,u_j,u_i\}=\frac{1}{2}u_i$,
\item[(SPG3)] $\{\widetilde{u}_i,\widetilde{u}_i,\widetilde{u}_j\}=\frac{1}{2}\widetilde{u}_j$,
$\{\widetilde{u}_j,\widetilde{u}_j,\widetilde{u}_i\}=\frac{1}{2}\widetilde{u}_i$,
\item[(SPG4)] $\{u_i,u_j,\widetilde{u}_i\}=-\frac{1}{2}\widetilde{u}_j$,
\item[(SPG5)] $\{u_j,\widetilde{u}_i,\widetilde{u}_j\}=-\frac{1}{2}u_i$,
\item[(SPG6)] All other products of elements from the spin grid are $0$.
\end{description}
In the case of finite odd dimensions (where $u_0$ is present) we
have, for all $i\in J$, the additional conditions (as exceptions of (SPG6))
\begin{description}
\item[(SPG7)] $\{u_0,u_0,u_i\}=u_i$, $\{u_i,u_i,u_0\}=\frac{1}{2}u_0$,
\item[(SPG8)] $\{u_0,u_0,\widetilde{u}_i\}=\widetilde{u}_i$, $\{\widetilde{u}_i,\widetilde{u}_i,u_0\}=\frac{1}{2}u_0$,
 \item[(SPG9)]$\{u_0,u_i,u_0\}=-\widetilde{u}_i,
 \{u_0,\widetilde{u}_i,u_0\}=-u_i.$
\end{description}
It is known (cf.\
\cite{dangfriedman-Classificationoftriplefactorsandapplications})
that every finite-dimensional spin factor is linearly spanned by a
spin grid (but not necessarily by a spin system).

Let $\fG:=\{u_i,\widetilde{u}_i:i\in I\}$ (resp.\
$\widetilde{\fG}:=\fG\cup\{u_0\}$) be a spin grid which spans the
$JC^*$-triple $Z$ and $1\in I$ an arbitrary index. If we define a
tripotent $v:=i(u_1+\widetilde{u}_1)$, Neal and Russo gave a method
how to construct from $\fG$ (resp.\ $\widetilde{\fG}$) and $v$ a
$JC^*$-triple system in
\cite{nealrusso-Contractiveprojectionsandoperatorspaces}, which is $JB^*$-triple isomorphic to $Z$ and
contains a spin system. First they have shown for the Peirce-2-space
$P_2^v(Z)$ of $v$ that $P_2^v(Z)=Z$ and that, if $\mathfrak{A}$ is
any von Neumann algebra containing $Z$, then
$P_2^v(\mathfrak{A})^{(v)}$ is a $C^*$-algebra TRO-isomorphic to
$P_2^v(\mathfrak{A})$ (the isomorphism is the identity mapping).
Moreover, they proved:

\begin{thm}[\cite{nealrusso-Contractiveprojectionsandoperatorspaces},
3.1]\label{AusGridMachSpinSystem} The space $P_2^v(Z)^{(v)}$ is the
linear span of a spin grid. More precisely, let
$s_j=u_j+\widetilde{u}_j, j\in I\setminus\{1\};$
$t_j:=i(u_j-\widetilde{u}_j), j\in I$. Then a spin system in the which linearly spans
$P_2^v(Z)^{(v)}$, is given by \[ \{s_j,t_k,v:j\in
I\setminus\{1\},k\in I\}
\]
or, if the spin factor is of odd finite dimension\[
\{s_j,t_k,v,u_0:j\in I\setminus\{1\},k\in I\}.
\]
\end{thm}

\begin{lemma}\label{PeirceraumVomTRO}
Let $T$ be a TRO and $v\in \Tri(T)$.
\begin{description}
\item[(a)] We have $P_2^v(T)=\{z\in
T:v\left(vz^*v\right)^*v=z\}.$
\item[(b)] Let $Z\subseteq B(H)$ be a $JC^*$-triple system and $T$ the TRO generated
by $Z$. If $Z=P_2^v(Z)$, then $T=P_2^v(T)$.
\item[(c)] If $v$ is a tripotent in the TRO $T$, then the
Peirce-2-space $P_2^v(T)$ is a subTRO of $T$.
\end{description}

\begin{proof}
\begin{description}
\item[(a)] Let $z\in T$ with $vv^*z+zv^*v=2z$. Then $vv^*$ and $v^*v$ are
projections with $vv^*zv^*v+zv^*v=2zv^*v$ and
$vv^*zv^*v+vv^*z=2vv^*z$. Thus we have $vv^*zv^*v=zv^*v=vv^*z$ and
therefore $vv^*zv^*v=\frac{1}{2}(vv^*z+zv^*v)=z$.

If $z\in Z$ with $vv^*zv^*v=z$, then $vv^*zv^*v=zv^*v$ and
$vv^*zv^*v=vv^*z$. We get $\frac{1}{2}(vv^*z+zv^*v)=vv^*zv^*v=z.$

\item[(b)] Let $x=z_1z_2^*z_3\ldots z_{2n} z_{2n+1}\in T$, with $z_j\in
Z=P_2^v(Z)$. By (a) we get $vv^*z_jv^*v=z_j$ and
$z_j=vv^*z_j=z_jvv^*$. Thus $vv^*xv^*v=\linebreak
(vv^*z_1)z_2^*z_3\ldots z_{2n}^*(z_{2n+1}v^*v)=z_1z_2^*z_3\ldots
z_{2n}^*z_{2n+1}=x$ and it follows that $x\in P^v_2(T)$.
\item[(c)] Let $a,b,c\in P_2^v(T)$, then
$vv^*ab^*cv^*v=vv^*a(vv^*bv^*v)^*cv^*v=\linebreak
(vv^*av^*v)b^*(vv^*cv^*v)=ab^*c$.
\end{description}
\end{proof}
\end{lemma}

As a first result we get an upper bound for the dimension of the
universal enveloping TRO of a spin system:

\begin{prop}\label{obereGrenzefuerspin}
Let $Z$ be a spin factor of dimension $k+1<\infty$. Then \[\dim
T^*(Z)\leq 2^k.\]
\begin{proof}
For $k=2n$ let
$\fG=\{u_1,\widetilde{u}_1,\ldots,u_n,\widetilde{u}_n\}$ (resp.\
$\fG=
\{u_1,\widetilde{u}_1,\ldots,u_n,\widetilde{u}_n\}\cup\{u_0\}$ for
$k=2n+1$) be a spin grid generating $Z$. Then $\rho_Z(\fG)$ is a
spin grid in $\rho_Z(Z)\subseteq T^*(Z)$. By Lemma
\ref{PeirceraumVomTRO} we have for $v:=i(u_1+\widetilde{u}_1)$ that
$P_2^v(T^*(Z))=T^*(Z)$, which is TRO-isomorphic to
$P_2^v(T^*(Z))^{(v)}$. The unital $C^*$-algebra
$P_2^v(T^*(Z))^{(v)}$ contains by Theorem
\ref{AusGridMachSpinSystem} a spin system $\{\id,s_1,\ldots,s_k\}$,
which generates it as a $C^*$-algebra. It is easy to see (cf.\
\cite{HancheOlsenStoermerJordanoperatoralgebras}, Remark 7.1.12) that
$P_2^v(T^*(Z))^{(v)}$ is linearly spanned by the $2^k$ elements
$s_{i_1}\ldots s_{i_j}$, where $1\leq i_1<i_2<\ldots<i_j$ and $0\leq
j\leq k$.
\end{proof}
\end{prop}

From the proof of Proposition \ref{obereGrenzefuerspin} we can deduce that the universal
enveloping TRO of a spin factor is TRO-isomorphic to its universal
enveloping $C^*$-algebra,
once we have shown that $\dim T^*(Z)=2^k$.

In Jordan-$C^*$-theory the following famous spin system appears\\
(cf.\ \cite{HancheOlsenStoermerJordanoperatoralgebras}, 6.2.1):

Let \[\sigma_1:=\left(
                  \begin{array}{cc}
                    1 & 0 \\
                    0 & -1 \\
                  \end{array}
                \right),\;\;\sigma_2:=\left(
                                        \begin{array}{cc}
                                          0 & 1 \\
                                          1 & 0 \\
                                        \end{array}
                                      \right)\;\text{ and }\;\sigma_3:=\left(
                                                                     \begin{array}{cc}
                                                                       0 & i \\
                                                                       -i & 0 \\
                                                                     \end{array}
                                                                   \right)
 \]
 be the Pauli spin matrices.

 For matrices $a=(\alpha_{i,j})\in \M_k$ and $b\in
 \M_l$ we define $a\otimes b:=(\alpha_{i,j}b)\in
 M_k(\M_l)=\M_{kl}$. If $a\in\M_k$ and $l\in\N$ we denote by $a^{\otimes l}$ the $l$-fold tensor product of $a$ with itself.

 The so called standard spin system, which is linearly generating a $\lr{k+1}$-dimensional spin factor in $\M_{2^n}$, when $k\leq2n$, is given via
 $\{\id,s_1,\ldots,s_k\}$ with
\begin{tabbing}
                $s_1:=\sigma_1\otimes\id^{\otimes (n-1)}$\;\;\;\;\;\;\;\;\;\;\;\;\;\;\;\;\;\;\;\;\;\;\;\;\;\;\;\;\;\;\;\;\;\;\;\;\;\;\=
                $s_2:=\sigma_2\otimes\id^{\otimes (n-1)},$\\\\
                $s_3:=\sigma_3\otimes\sigma_1\otimes\id^{\otimes (n-2)},$
                \> $s_4:=\sigma_3\otimes\sigma_2\otimes\id^{\otimes (n-2)},$\\\\
                $s_{2l+1}:=\sigma_3^{\otimes l}\otimes\;\sigma_1\otimes\id^{\otimes (n-l-1)},$
                \>$s_{2l+2}:=\sigma_3^{\otimes l}\otimes\;\sigma_1\otimes\id^{\otimes (n-l-1)} $

             \end{tabbing}
             for $1\leq l\leq n-1$.
\begin{lemma}
Let $S=\{\id,s_1,\ldots,s_k\}$ be the standard spin system. If
$k=2n$, then the TRO generated by $S$ in $\M_{2^n}$ is
$\M_{2^n}$. If $k=2n-1$ then the generated TRO is
TRO-isomorphic to $\M_{2^{n-1}}\oplus
\M_{2^{n-1}}$.
\begin{proof}
Let $T$ be the TRO generated by $S$.

Let $k=2n$. It suffices to show that the $3k$ elements
\[
  a_j:=\id^{\otimes (j-1)}\otimes\sigma_1\otimes\id^{\otimes (n-j)},\;\;\;\;
   b_j:=\id^{\otimes (j-1)}\otimes\sigma_2\otimes\id^{\otimes (n-j)},\]
  \[ c_j:=\id^{\otimes (j-1)}\otimes\sigma_1\otimes\id^{\otimes (n-j)}
\]
for every $j=1,\ldots,k$ are elements of $T$, since $a_j,b_j,c_j$
and $\id\otimes\ldots\otimes\id$ span
$\mathbb{C}\otimes\ldots\otimes\mathbb{C}\otimes
\M_2\otimes\mathbb{C}\otimes\ldots\otimes\mathbb{C}$.

Obviously $a_1=s_1\in T$. Suppose we have shown $a_j\in T$ for a fixed
$j\geq 1$, then
\begin{align*}
s_{2j}s_{2j+1}^*a_j &=\id^{\otimes (j-1)}\otimes\sigma_2\sigma_3\sigma_1\otimes\sigma_1\id^{\otimes (n-j-1)}\\
&=ia_{j+1}.
\end{align*}
Similarly we have $b_1=s_2\in T$. If we have shown for a fixed
$j\geq 1$ that $b_j\in T$, then
\[
s_{2j}s_{2j+2}^*a_j =ib_{j+1}.
\]
Another easy induction shows that $c_j\in T$ for all $j=1,\ldots,n$.

If $k=2n-1$ we have $a_n\in T$, $b_n,c_n\notin T$. Since $\sigma_1$
and $\id\otimes\ldots\otimes\id$ generate the diagonal matrices, the
statement is clear.

Alternatively we could argue that $T$ contains the identity
so $T$ has to be a $C^*$-algebra. Then the statement follows from
\cite{HancheOlsenStoermerJordanoperatoralgebras}, Theorem 6.2.2.
\end{proof}
\end{lemma}

\begin{thm}\label{enveloping TRO for Spin}
For the universal enveloping TRO of a spin factor $Z$ with $\dim
Z=k+1$ we have
\[T^*(Z)=\begin{cases}

  \M_{2^{n-1}}\oplus \M_{2^{n-1}}  & \text{if }k=2n-1,\\
  \M_{2^{n}} & \text{if }k=2n.

\end{cases}\]

\begin{proof}
The $JC^*$-triple system $Z$ is $JB^*$-isomorphic to the
$JC^*$-algebra $J$ linearly generated by the standard spin system
$\{1,s_1,\ldots,s_k\}$.
 By the universal property of $T^*(Z)$ we get, since $J$ generates
 $\M_{2^{n-1}}\oplus \M_{2^{n-1}} \text{ if
 }k=2n-1$ (respectively $\M_{2^{n}}  \text{ if }k=2n$) as
 a TRO, a surjective TRO-homomorphism from $T^*(Z)$ onto $\M_{2^{n-1}}\oplus\M_{2^{n-1}}$ if $k=2n-1$ (respectively $\M_{2^{n}}$ if $k=2n$).
 By Proposition
 \ref{obereGrenzefuerspin} this has to be an isomorphism.
\end{proof}

\end{thm}

\subsection{Factors of type III}\label{Factors of type III}
A hermitian grid is a family $\{u_{ij}:i,j\in I\}$ of tripotents in
$Z$ such that for all $i,j,k,l\in I$:
\begin{description}
\item[(HG1)] $u_{ij}=u_{ji}$ for all $i,j\in I$.
\item[(HG2)] $\{u_{kl},u_{kl},u_{ij}\}=0$ if
$\{i,j\}\cap\{k,l\}=\emptyset$.
\item[(HG3)] $\{u_{ii},u_{ii},u_{ij}\}=\frac{1}{2}u_{ij}$,
$\{u_{ij},u_{ij},u_{ii}\}=u_{ii}$ if $i\neq j$.
\item[(HG4)] $\{u_{ij},u_{ij},u_{jk}\}=\frac{1}{2}u_{jk}$,
$\{u_{jk},u_{jk},u_{ij}\}=\frac{1}{2}u_{ij}$ if $i,j,k$ are pairwise
distinct.
\item[(HG5)] $\{u_{ij},u_{jk},u_{kl}\}=\frac{1}{2}u_{il}$ if $i\neq l$.
\item[(HG6)] $\{u_{ij},u_{jk},u_{ki}\}=u_{ii}$ if at least two of
these tripotents are distinct.
\item[(HG7)] All other products of elements from the hermitian grid are $0$.
\end{description}
Let $Z$ be a finite-dimensional TRO. Then the direct sum \[T=\bigoplus_{\alpha=1}^r\M_{n_\alpha,m_\alpha}\]
can be described by so called \textbf{rectangular matrix units}\index{rectangular matrix units}: Let $E(\alpha,i,j):=E_{i,j}\in \M_{n_\alpha,m_\alpha}$ be the matrix in $\M_{n_\alpha,m_\alpha}$ which is $0$ everywhere except $1$ in the $(i,j)$-component for all $1\leq i\leq n_\alpha$, $1\leq j\leq m_\alpha$ and $\alpha\in\{1,\ldots,r\}$.
Put $$
e^{(\alpha)}_{i,j}:=(0,\ldots,0,E(\alpha,i,j),0,\ldots,0)\in T,
$$where $E(\alpha,i,j)$ is in the $\alpha$th summand. The rectangular matrix units satisfy
\begin{description}
  \item[(i)] $e_{i,j}^{(\alpha)}\left(e_{l,j}^{(\alpha)}\right)^*e_{l,k}^{(\alpha)}=e_{i,k}^{(\alpha)}$.
  \item[(ii)] $e_{i,j}^{(\alpha)}\left(e_{n,m}^{(\beta)}\right)^*e_{p,q}^{(\gamma)}=0$ for $j\neq m$, $n\neq p$, $\alpha\neq\beta$ or $\beta\neq \gamma$.
  \item[(iii)] $T=\lin\{e_{i,j}^{(\alpha)}:1\leq\alpha\leq r,\; 1\leq i\leq n_\alpha,\; 1\leq j\leq m_\alpha\}$.
\end{description}
If $U$ is another TRO which contains elements $f_{i,j}^{(\beta)}$ satisfying the analogues of $(i)$--$(iii)$ for $1\leq i\leq n_\alpha$, $1\leq j\leq m_\alpha$ and $\alpha,\beta\in\{1,\ldots,r\}$, then it is easy to see that the mapping sending $e_{i,j}^{(\alpha)}$ to $f_{i,j}^{(\alpha)}$ for $1\leq i\leq n_\alpha$, $1\leq j\leq m_\alpha$ and $\alpha\in\{1,\ldots,r\}$ is a TRO-isomorphism.

Let $Z$ be a finite-dimensional $JC^*$-triple system spanned by a
hermitian grid $\{u_{ij}:1\leq i,j\leq n\}$ and $T$ the TRO
generated by this grid. Define $$e_{ij}:=u_{ii}(\sum_{k=1}^n
u_{kk})^*u_{ji}\in T$$ for $1\leq i,j\leq n,$. From
\cite{nealrusso-Contractiveprojectionsandoperatorspaces}, Lemma 3.2 (a) we can conclude
 that $\{e_{ij}\}$
forms a system of rectangular matrix units in $T$. We get that $$T^*(Z)\simeq \M_n.$$

\subsection{Factors of type II}\label{Factors of type II}

 A symplectic grid is a family $\{u_{ij}:i,j\in I,
i\neq j\}$ of minimal tripotents such that for all $i,j,k,l\in I$
\begin{description}
\item[(SYG1)] $u_{ij}=-u_{ji}$ for $i\neq j$.
\item[(SYG2)] $\{u_{ij},u_{ij},u_{kl}\}=\frac{1}{2}u_{kl}$,
$\{u_{kl},u_{kl},u_{ij}\}=\frac{1}{2}u_{ij}$ for
$\{i,j\}\cap\{k,l\}\neq\emptyset$.
\item[(SYG3)] $\{u_{kl},u_{kl},u_{ij}\}=0$ if
$\{i,j\}\cap\{k,l\}=\emptyset$.
\item[(SYG4)] $\{u_{ij},u_{il},u_{kl}\}=\frac{1}{2}u_{kj}$ for
$i,j,k,l$ pairwise distinct.
\item[(SYG5)] All other triple products in the symplectic grid are
$0$.
\end{description}
The standard example of a finite-dimensional symplectic grid is the
collection $\{U_{i,j}:1\leq i,j\leq n, i\neq j\}\subseteq
\M_n$, where $U_{i,j}$, for $i<j$, is a complex $n\times
n$- matrix, which is 0 everywhere except for the $(i,j)$-entry,
which is 1 and the $(j,i)$-entry, which is $-1$. This grid spans
linearly the $JC^*$-triple system $\{A\in \M_n:A^t=-A\}$
of skew-symmetric $n \times n$ matrices; its TRO span is
$\M_n$.

Let $\fG:=\{u_{ij}:i,j\in I, i\neq j\}$ be a symplectic grid, $Z$
the $JC^*$-triple system spanned by $\fG$ and $T$ the TRO generated
by it. Since for $\dim Z=3$ $Z$ is $JB^*$-triple isomorphic to a
type I Cartan factor and for $\dim Z=6$ it is $JB^*$-triple
isomorphic to a type IV Cartan factor, both covered in other
sections, let $\dim Z\geq 10$.

If we define
$$e_{ii}:=u_{ik}u_{kl}^*u_{il}$$ and
$$e_{ij}:=e_{ii}e_{ii}^*u_{ij}e_{jj}^*e_{jj},$$
for
$1\leq i,j,k,l \leq n$ pairwise distinct,
we get with \cite{nealrusso-Contractiveprojectionsandoperatorspaces},
Lemma 4.1 and Lemma 4.3 that
the elements $e_{ii}$ and $e_{ij}$ are well-defined and that for $v:= \sum e_{kk}$ we have $$ve_{ij}^*v=e_{ji}\;\;\;\;\;\;\;\;\;\text{ and }
\;\;\;\;\;\;\;\;\;e_{ij}v^*e_{kl}=\delta_{jk}e_{il}.$$
Using this we see that
\begin{align*}
  e_{ij}e_{kl}^*e_{mn} 
   &=e_{ij}v^*e_{lk}v^*e_{mn}\displaybreak[0] \\
   &=\delta_{jl}\delta_{km} e_{in},
\end{align*}
which shows that $\{e_{ij}\}$ is a set of rectangular matrix units.
\begin{thm}
If $Z$ is a $JC^*$-triple system spanned by a symplectic grid
with $\dim Z\geq 10$,
then
\[T^*(Z)=\M_n.\]
\end{thm}

\subsection{Factors of type I}\label{Factors of type I}
Let $\Delta$ and $\Sigma$ be two index sets.
 A rectangular grid is a family $\{u_{ij}:i\in \Delta,
j\in \Sigma\}$ of minimal tripotents such that
\begin{description}
\item[(RG1)]$\{u_{il},u_{il},u_{jk}\}=0$ if
 $i\neq j, k\neq l$.
\item[(RG2)] $\{u_{il},u_{il},u_{jk}\}=\frac{1}{2}u_{jk}$,
$\{u_{jk},u_{jk},u_{il}\}=\frac{1}{2}u_{il}$ if either $j=i,k\neq l$ or $j\neq
i,k=l$.
\item[(RG3)] $\{u_{jk},u_{jl},u_{il}\}=\frac{1}{2}u_{ik}$ if $j\neq
i$ and $k\neq l$.
\item[(RG4)] All other triple products in the rectangular grid equal
$0$.
\end{description}

Let $Z$ be the $JC^*$-triple system generated by a finite
rectangular grid. We assume that $Z$ is finite-dimensional and
hence $JB^*$-triple isomorphic to $\M_{n,m}$ with
$m=|\Delta|$ and $n=|\Sigma|$.

We first exclude some candidates for $T^*(Z)$:
\begin{lemma}\label{TROHuelleNichtIsomZuMatrizen}
For the $JC^*$-triple system $Z=\M_{n,m}$ its universal
enveloping TRO $T^*(Z)$ is neither TRO-isomorphic to
$\M_{n,m}$ nor to $\M_{m,n}$.
\begin{proof}
Assume that $T^*(Z)$ is TRO-isomorphic to $\M_{n,m}$. Let
$^t:\M_{n,m}\to \M_{m,n}$ be the transposition
mapping. According to the universal property of $T^*(Z)$ there is a
mapping $T^*(^t)$ such that
\[
\xymatrix{&  \M_{n,m}\ar[rd]^{T^*(^t)}
\\
\M_{n,m}\ar[ru]^-{\rho_Z}\ar[rr]^{^t}&&
\M_{m,n} }
\] commutes. Since $\rho_Z$ is bijective there is a TRO-isomorphism
$T^*(\rho_Z):\M_{n,m}\to \M_{n,m}$ with
$T^*(\rho_Z)\circ\rho_Z=\id$. This means $T^*(\rho_Z)=\rho_Z^{-1}$,
in particular $\rho_Z$ is a complete isometry.

Since $\rho_Z$ and $^t$ are bijective the same holds for $T^*(^t)$
and it follows that $^t$ is a complete isometry. We get a contradiction because $^t$ is not
even completely bounded. The other statement
can be proved analogously.
\end{proof}
\end{lemma}

\begin{lemma}[\cite{nealrusso-Contractiveprojectionsandoperatorspaces}, Lemma 5.1 (b), Lemma 5.2
(b)]\label{isomorphzuTRO} Let $\{u_{ij}\}$ be a rectangular grid
spanning $Z$.
\begin{description}
  \item[(a)] If for $i\in\Delta,$ $k,l\in \Sigma$, where $k\neq l$, we
  have $u_{il}u_{ik}^*=0$  or for $i,j\in\Delta,$ $k\in \Sigma$, where $i\neq
  j$, we have $u_{ik}^*u_{jk}=0$, then $Z$ is TRO-isomorphic to
  $\M_{n,m}$.
  \item[(b)] If for $i\in\Delta,$ $k,l\in \Sigma$, where $k\neq l$, we
  have $u_{il}^*u_{ik}=0$  or for $i,j\in\Delta,$ $k\in \Sigma$, where $i\neq
  j$, we have $u_{ik}u_{jk}^*=0$, then $Z$ is TRO-isomorphic to
  $\M_{m,n}$.
\end{description}
\end{lemma}

By this we get
\begin{lemma}
Let $\{e_{ij}\}$ be a rectangular grid spanning $\rho_Z(Z)\subseteq
T^*(Z)$, then we have
\begin{equation}\label{notwFuerRecGrid1}
    e_{ik}e_{il}^*\neq0\text{ and }e_{ik}^*e_{il}\neq
    0\text{ for all }i\in\Delta,k,l\in\Sigma
\end{equation}
as well as
\begin{equation}\label{notwFuerRecGrid2}
    e_{ik}e_{jk}^*\neq0\text{ and }e_{ik}^*e_{jk}\neq
    0\text{ for all }i,j\in\Delta,k\in\Sigma.
\end{equation}
\begin{proof}
If one of these conditions is not fulfilled we get by Lemma
\ref{isomorphzuTRO} and since $\rho_Z(Z)$ generates $T^*(Z)$ as a
TRO, that $\rho_Z(Z)=T^*(Z)$ and hence is isomorphic to
$\M_{n,m}$ respectively $\M_{m,n}$. But this
is a contradiction to Lemma \ref{TROHuelleNichtIsomZuMatrizen}.
\end{proof}
\end{lemma}

\begin{lemma}
Let $\rank Z\geq 2$ and $\{e_{ij}\}$ be a rectangular grid spanning
$\rho_Z(Z)$, then
\[p:=\sum_{i\in\Delta}\prod_{j\in\Sigma}e_{ij}e_{ij}^*\in C^*(Z)\]
is a sum of non-zero orthogonal projections. We have:

\[pT^*(Z)\subseteq T^*(Z),\;\; (1-p)T^*(Z)\subseteq T^*(Z).\]
\begin{proof}Since (\ref{notwFuerRecGrid1}) and
(\ref{notwFuerRecGrid2}) hold we can use
\cite{nealrusso-Contractiveprojectionsandoperatorspaces}, Lemma 5.5
and get that $\prod_{j\in\Sigma}e_{ij}e_{ij}^*\neq 0$ are orthogonal
projections for all $i\in \Delta$.

The fact that $p$ leaves $T^*(Z)$ invariant is obvious. 
\end{proof}
\end{lemma}

\begin{lemma}\label{ProjektionsTROs}
For all $i,k,a\in \Delta$, $j,l,b\in\Sigma$ we have \[
pe_{ij}(pe_{kl})^*pe_{ab}=pe_{ij}e_{kl}^*pe_{ab}\in\lin\{pe_{ij}\}
\]
and for $q:=(1-p)$
\[
qe_{ij}(qe_{kl})^*qe_{ab}=qe_{ij}e_{kl}^*qe_{ab}\in\lin\{qe_{ij}\}.
\]
\begin{proof}
Since $\{e_{ij}\}$ is a rectangular grid we know for $i\neq k$ and
$j\neq l$ that
\begin{equation*}\label{OrthogImRrcGrid}
    e_{ij}e_{kl}^*=0\;\;\;\;\;\;\text{ and }\;\;\;\;\;\;e_{ij}^*e_{kl}=0
\end{equation*}
and therefore, for $i\neq k$ and $j\neq l$,
\begin{equation}\label{ProjGrid1}
  pe_{il}(pe_{kl})^* = pe_{il}e_{kl}^*p = 0\end{equation}
as well as
\begin{align}\label{ProjGrid2}
\nonumber  (pe_{il})^*pe_{kl} &= e_{il}^*pe_{kl} \displaybreak[0]\\
\nonumber   &= e_{il}^*\lr{\sum_{\alpha\in\Delta}\prod_{\beta\in\Sigma}e_{\alpha \beta}e_{\alpha \beta}^*}e_{kl}\displaybreak[0] \\
\nonumber   &= e_{il}^*e_{i1}e_{i1}^*\ldots
\nonumber   e_{in}\underbrace{e_{in}^*e_{kl}}_{=0\text{ if } n\neq l}\displaybreak[0]\\
   &= 0,
\end{align}
since the range projections of collinear tripotents commute by
\cite{nealrusso-Contractiveprojectionsandoperatorspaces}, Lem\-ma 5.4.

Equation (\ref{ProjGrid1}) and (\ref{ProjGrid2}) lead us to the fact that we only
have to prove (for arbitrary $a\in \Delta,b\in\Sigma$) that
\begin{tabbing}
  $\bullet$ $pe_{ik}(pe_{il})^*pe_{ab}$ \;\;\;\;\;\=$k\neq l$\;\;\;\;\;\;\;\;\;\;\;\;\;\;\=
  $\bullet$ $pe_{jk}(pe_{ik})^*pe_{ab}$ \;\;\;\;\;\=$i\neq j$\\\\
  $\bullet$ $pe_{il}(pe_{il})^*pe_{ab}$
  \>\>$\bullet$ $pe_{ab}(pe_{il})^*pe_{ik}$ \>$k\neq l$\\\\
  $\bullet$ $pe_{ab}(pe_{il})^*pe_{jl}$ \>$i\neq j$ \>
  $\bullet$ $pe_{ab}(pe_{il})^*pe_{il}$\\
\end{tabbing}
are elements of $\lin\{pe_{ij}\}$.

Using (\ref{ProjGrid1}) and (\ref{ProjGrid2}) again, we have to prove
this in the following cases:

\begin{tabbing}
$\bullet$ $pe_{ik}(pe_{il})^*pe_{ib}$ \;\;\;\;\;\=$k\neq l$, $k\neq
b\neq
 l$\;\;\;\;\;\;\;\=
$\bullet$ $pe_{ik}(pe_{il})^*pe_{ik}$ \;\;\;\;\;\=$k\neq l$\\\\
$\bullet$ $pe_{ik}(pe_{il})^*pe_{il}$ \>$k\neq l$
\>$\bullet$ $pe_{ik}(pe_{il})^*pe_{al}$ \>$k\neq l$, $a\neq i$\\\\
$\bullet$ $pe_{jk}(pe_{ik})^*pe_{ib}$ \>$b\neq k$, $i\neq j$
\>$\bullet$ $pe_{jk}(pe_{ik})^*pe_{ik}$ \>$i\neq j$\\\\
$\bullet$ $pe_{jk}(pe_{ik})^*pe_{jk}$ \>$i\neq j$ \>$\bullet$
$pe_{jk}(pe_{ik})^*pe_{ak}$ \>$i\neq j$, $a\neq
i$\\\\
$\bullet$ $pe_{il}(pe_{il})^*pe_{ib}$ \>$b\neq l$
\>$\bullet$ $pe_{il}(pe_{il})^*pe_{al}$ \>$a\neq i$\\\\
$\bullet$ $pe_{il}(pe_{il})^*pe_{il}$.
\end{tabbing}
 We obtain a similar list for $q$. Luckily, Neal and Russo calculated all these products to show that $\{pe_{ij}\}$ is a
 rectangular grid (cf.\ the proof of \cite{nealrusso-Contractiveprojectionsandoperatorspaces}, Lemma
 5.6) and it is true that all of them are elements of
 $\{pe_{ij}\}$. One can show by similar methods that all products in
 the list for $q$ are elements of the rectangular grid
 $\{(1-p)e_{ij}\}$.
\end{proof}
\end{lemma}

\begin{prop}\label{ObereGrenzeTyp1}
If $\rank Z\geq 2$ we have for the universal enveloping TRO of $Z$
\[ T^*(Z)=\lin \{pe_{ij},(1-p)e_{ij}:1\leq i\leq n, 1\leq i\leq m\}
\]especially \[\dim T^*(Z)\leq 2nm.\]
\begin{proof}
The rectangular grid $\{e_{ij}\}$ spans $\rho_{Z}(Z)$ which
generates $T^*(Z)$ as a TRO, so an element $x\in T^*(Z)$ has to be
of the form
\[x=\sum_{\alpha=1}^n\lambda_\alpha e_1^\alpha(e_2^\alpha)^*e_3^\alpha\ldots
(e_{2n}^\alpha)^*e_{2k_\alpha+1}^\alpha,\] with $e_1^\alpha,\ldots,
e_{2k_{\alpha}+1}^\alpha\in\{e_{ij}\}$, $\lambda_{\alpha}\in
\mathbb{C}$ and $k_{\alpha}\in \mathbb{N}$ for all $1\leq \alpha\leq
n$, $n\in\mathbb{N}$. Let $e_1,\ldots ,e_{2n+1}$ and $e:=e_1e_2^*e_3\ldots e_{2n}e_{2n+1}^*\in
T^*(Z)$, then
\begin{align*}
  e &= \lr{pe_1+(1-p)e_1}\lr{pe_2+(1-p)e_2}^*\ldots\lr{pe_{2n+1}+(1-p)e_{2n+1}}\displaybreak[0] \\
   &= pe_1\lr{pe_2}^*\ldots pe_{2n+1}
   +(1-p)e_1\lr{(1-p)e_2}^*\ldots(1-p)e_{2n+1}\displaybreak[0]\\
   &+\text{ mixed terms in $p$ and $(1-p)$}\displaybreak[0]\\
   &= pe_1\lr{pe_2}^*\ldots pe_{2n+1}
   +(1-p)e_1\lr{(1-p)e_2}^*\ldots(1-p)e_{2n+1},
\end{align*}
since $\{pe_{ij}\}\perp\{(1-p)e_{ij}\}$ by Lemma
\cite{nealrusso-Contractiveprojectionsandoperatorspaces}, Lemma 5.6.
An inductive use of Lemma \ref{ProjektionsTROs} gives us $e\in
\{pe_{ij},(1-p)e_{ij}:1\leq i\leq n, 1\leq i\leq m\}$.
\end{proof}
\end{prop}

\begin{thm}
Let $Z$ be a $JC^*$-triple system of $\rank\geq 2$ and isomorphic
to a finite-dimensional Cartan factor of type I. Let $\{u_{ij};1\leq
i\leq n,1\leq j\leq m\}$ be a grid spanning $Z$. Then
\[T^*(Z)=\M_{n,m}\oplus \M_{m,n}.\]
\begin{proof}
We identify $Z$ with $\M_{n,m}$. The mapping
$\Phi:\M_{n,m}\to\M_{n,m}\oplus \M_{m,n}$, $A\mapsto
\left(A,A^t\right)$ is a $JB^*$-triple isomorphism onto a $JB^*$-subtriple of
$\M_{n,m}\oplus \M_{m,n}$ which generates
$\M_{n,m}\oplus \M_{m,n}$ as a TRO. Since by
\ref{ObereGrenzeTyp1} $\dim T^*(Z)\leq 2nm$ the induced mapping
$T^*(\Phi):T^*(Z)\to \M_{n,m}\oplus \M_{m,n}$
has to be a TRO isomorphism.
\end{proof}
\end{thm}

For the rest of this section we assume that $\rank Z=1$ and $Z$ is
of finite dimensions. This implies, that if $\{u_{ij};1\leq i\leq
n,1\leq j\leq m\}$ is a rectangular grid spanning $Z$ then $n$ or
$m$ have to be equal to $1$. In this special case the definition of
a rectangular grid becomes simpler:

A finite rectangular grid of rank $1$ is a set $\{u_1,\ldots,u_n\}$
of tripotents with
\begin{description}
\item[(RG'1)] $\{u_i,u_j,u_i\}=0$ for $i\neq j$.
\item[(RG'2)] $\{u_i,u_i,u_k\}=\fr{2}u_k$ for $i\neq k$.
\item[(RG'3)] All other products are $0$.
\end{description}

Let $Z$ be a $n$-dimensional type $1$ Cartan factor of $\rank 1$. We
fix a finite rectangular grid $\{e_1,\ldots,e_n\}$
of rank $1$ spanning $\rho_Z(Z)\subseteq T^*(Z)$.

\begin{lemma}\label{obereGrenzefuerRang1}
Let $Z$ be as above, then
\[\dim T^*(Z)\leq \sum_{k=1}^n\begin{pmatrix}n\\ k-1\end{pmatrix}\begin{pmatrix}n\\ k\end{pmatrix}.\]
\begin{proof}
Using the grid properties (RG'1),(RG'2),(RG'3) we show that
\begin{align*}
T^*(Z)=\lin\{&e_{i_1}e_{i_2}^*e_{i_3}\ldots
e_{i_{2k}}^*e_{i_{2k+1}}:
i_j<i_{j+2},1\leq j\leq 2k-1,\displaybreak[0]\\
&  0\leq k\leq \frac{1}{2}(n-1)\}.
\end{align*}
 For a fixed $k$ we have $\begin{pmatrix}n\\
k-1\end{pmatrix}\begin{pmatrix}n\\ k\end{pmatrix}$ choices for
$e_{i_1}e_{i_2}^*e_{i_3}\ldots e_{i_{2k}}^*e_{i_{2k+1}}$. This is
true because
$i_j<i_{j+2}$. We have $\begin{pmatrix}n\\
k\end{pmatrix}$ choices for $i_1< i_3<\ldots<i_{2k+1}$ and $\begin{pmatrix}n\\
k-1\end{pmatrix}$ choices for $i_2< i_4<\ldots<i_{2k}$.

To prove that $T^*(Z)$ is the above mentioned linear span we give an induction which takes $x=e_{i_1}e_{i_2}^*e_{i_3}\ldots
e_{i_{2k}}^*e_{i_{2k+1}}\in T^*(Z)$ and rearranges the grid elements
such that $x$ is a sum of elements of the form
$e_{j_1}e_{j_2}^*e_{j_3}\ldots e_{j_{2k}}^*e_{j_{2k+1}}$ with
$j_1\leq j_3\leq\ldots\leq j_{2l+1}$ and $j_2\leq j_4\leq\ldots\leq
j_{2l}$. Since the grid elements are tripotents we can assume that
we do not have three equal indices in a row. If we have the
case $e_{\alpha}e_\beta^*e_\alpha$, where $\alpha\neq\beta$
this equals $0$ by the minimality of the tripotents (cf.\ (RG'1)).
Therefore $j_a<j_{a+2}$ for all $1\leq a\leq 2l-1$. Especially
$l\leq \frac{1}{2}(n-1)$.

So let $x=e_{i_1}e_{i_2}^*e_{i_3}\ldots e_{i_{2k}}^*e_{i_{2k+1}}\in
T^*(Z)$. Since the $e_{i_a}$ are all minimal tripotents we can
assume $e_{i_a}\neq e_{i_{a+2}}$.

For $k=0$ nothing is to prove. Additionally we prove the case when
$k=1$. Let $x=e_{i_1}e_{i_2}^*e_{i_3}$.
\begin{description}
\item If $i_1<i_3$ we are done.

\item If $i_1=i_2>i_3$ we can use (RG'2) and get
$x=e_{i_3}-e_{i_3}e_{i_1}^*e_{i_1}.$


\item If $i_1>i_2=i_3$ we can also use (RG'2) and get
$x=e_{i_1}-e_{i_2}e_{i_2}^*e_{i_1}$.
\item If $i_1\neq i_2\neq i_3$:
\begin{description}
\item If $i_1> i_3$ we can use (RG'3) and we deduce $x=-e_{i_3}e_{i_2}^*e_{i_1}$.
\end{description}
\end{description}
Now we assume that we have shown the statement for $2k+1\in
\mathbb{N}$, $2k+3\leq n$ and for all lesser indices. If we apply
our induction statement to the first $2k+1$ grid elements in the
product and then apply the beginning of the induction to all the
last three elements of the products in the resulting sum, then one can easily convince himself that in at most three repetitions of
this procedure we get the desired  form for $x$.
\end{proof}
\end{lemma}

Again we have to give a faithful representation $T$ of $T^*(Z)$. This happens
to be more complicated than in the other cases. Again we can use the work of Neal and Russo. In
\cite{nealrusso-Contractiveprojectionsandoperatorspaces} they showed
that a $JC^*$-triple system, which is linearly spanned by a finite
rectangular grid of rank $1$ with $n$ elements, has to be completely
isometric (especially $JB^*$-triple isomorphic) to one of the spaces
$H^k_n$, where $k=1,\ldots,n$, that are generalizations of the row
and column Hilbert space.

We recall the construction of the spaces $H^k_n$ (cf.\
\cite{nealrusso-Contractiveprojectionsandoperatorspaces}, Section 6
and 7 or
\cite{NealRusso-RepresentationofcontractivelycomplementedHilbertianoperatorspacesontheFockspace},
Section 1). Let $1\leq k\leq n$ and $I,J$ be subsets of
$\{1,\ldots,n\}$ such that $I$ has $k-1$ and $J$ has $n-k$ elements.
There are $q_k:=\begin{pmatrix}n\\ k-1\end{pmatrix}$ choices for $I$
and $p_k:=\begin{pmatrix}n\\ n-k\end{pmatrix}=\begin{pmatrix}n\\
k\end{pmatrix}$ choices for $J$. We assume that the collections
$\mathcal{I}:=\{I_1,\ldots,I_{q_k}\}$ and
$\mathcal{J}:=\{J_1,\ldots,J_{p_k}\}$ of such sets are ordered
lexicographically. Let $e_{I_1},\ldots,e_{I_{q_k}}$ and
$e_{J_1},\ldots,e_{J_{p_k}}$ be the canonical bases of
$\mathbb{C}^{p_k}$ and $\mathbb{C}^{q_k}$. We can define an element
in $\M_{p_k,q_k}$ by $E_{I,J}:=E_{i,j}$, when
$I=I_i\in\mathcal{I}$ and $J=J_j\in\mathcal{J}$. The space $H^k_n$
is the linear span of matrices $b^{n,k}_i$, where $1 \leq i \leq n$,
given by
\begin{equation}\label{H^k_n}
b^{n,k}_i:=\sum_{I\cap J=\emptyset, (I\cup J)^c=\{i\}}
\sgn(I,i,J)E_{J,I},
\end{equation}
where $\sgn(I,i,J)$ is the signature of the permutation
taking\linebreak
 $(i_1,\ldots,i_{k-1},i,j_1,\ldots,j_{n-k})$ to
$(1,\ldots,n)$, when $I=\{i_1,\ldots,i_{k-1}\}$, where
$i_1<i_2<\ldots<i_{k-1}$, and $J=\{j_1,\ldots,j_{n-k}\}$ and where
$j_1<j_2<\ldots<j_{n-k}$.

One can show that the TRO spanned by $b^{n,k}_1,\ldots,b^{n,k}_n$
equals $\M_{p_k,q_k}$, so if we represent our
$JC^*$-triple system $Z$ as $\bigoplus_{k=1}^nH^k_n$ we get with
Lemma \ref{obereGrenzefuerRang1}:
\begin{thm}
If $Z$ is a $JC^*$-triple system spanned by a finite rectangular
grid of rank $1$, then
\[T^*(Z)=\bigoplus_{k=1}^n\M_{p_k,q_k},\]
where $p_k=\begin{pmatrix}n\\
k\end{pmatrix}$ and $q_k=\begin{pmatrix}n\\ k-1\end{pmatrix}$ for all
$k=1,\ldots,n$.
\end{thm}
With this result the list of universal enveloping TROs of the finite-dimensional Cartan factors is complete.

\section{The Radical}

We use the theory of reversibility developed in \cite{BunceFeelyTimoneyI} to prove some facts for the universal enveloping TRO of a universally reversible TRO $T$. We consider the case in which a universally reversible TRO $T$ contains an ideal of codimension $1$ which is not covered in \cite{BunceFeelyTimoneyI}. We show that there exists an ideal $\ra(T)$ in $T$ which is universally reversible and which does not contain an ideal of codimension $1$ itself, such that $T/\ra(T)$ is an Abelian $JB^*$-triple system. We obtain an exact sequence
\begin{align*}
0\longrightarrow \ra(T)\oplus\theta(\ra(T))\longrightarrow  T^*(T)\longrightarrow C_0^{\mathbb{T}}(\Epi(T/\ra(T),\C))\longrightarrow 0,
\end{align*}
where the notation is given below.

We adopt the following definition from \cite{BunceFeelyTimoneyI}. It is the generalization of reversibility of $JC^*$-algebras.

\begin{defi}
A $JC^*$-triple system $Z\sub B(H)$ is said to be \textbf{reversible}\index{reversible} if \[\fr{2}\lr{x_1x_2^*x_3\ldots x_{2n}^*x_{2n+1}+x_{2n+1}x_{2n}^*\ldots x_3x_2^*x_1}\in Z\]
for all $x_1,\ldots,x_n\in Z$ and $n\in \N$. We call a $JC^*$-triple system \textbf{universally reversible}\index{universally reversible} if it is reversible in every representation.
\end{defi}
Obviously every TRO, and therefore every $C^*$-algebra, is reversible (but not necessarily universally reversible, since we have to cope with $JB^*$-triple homomorphisms).
A $JC^*$-triple system is universally reversible if and only if it is reversible when embedded in its universal enveloping TRO.

\begin{lemma}[\cite{BunceFeelyTimoneyI}, Theorem 4.4]\label{Tro identifizierer}
Let $Z$ be a universally reversible $JC^*$-triple system and let $\varphi:Z\to B(H)$ be an injective triple homomorphism. Suppose there exists a TRO antiautomorphism $\Psi$ of the TRO-span $\TRO(\varphi(Z))$ such that $\Psi\circ\varphi=\varphi$, then $T^*(\varphi):T^*(Z)\to \TRO(\varphi(Z))$ is a TRO-isomorphism.
\end{lemma}

\begin{lemma}[\cite{BunceFeelyTimoneyI}, Corollary 4.5]\label{codimension Tro}
Let $T$ be a universally reversible TRO in a $C^*$-algebra $\fA$. Suppose $T$ has no TRO-ideals of codimension $1$ and there is a TRO antiautomorphism $\theta:\fA\to\fA$ of order $2$. Then $T^*(T)\simeq T\oplus \theta(T)$ with universal embedding $a\mapsto (a,\theta(a))$.
\end{lemma}

In order to establish the announced generalization of Lemma \ref{codimension Tro} we define an ideal such that the quotient of $T$ by this ideal is Abelian. We are first recalling some known facts about Abelian $JB^*$-triple systems which allow us to compute the universal enveloping TRO of a general Abelian triple. Afterwards we show that every ideal of a universal reversible $JC^*$-triple system is universally reversible.

Recall that a $JB^*$-triple system $Z$ is called Abelian, if
\begin{equation*}
    \{\{a,b,c\},d,e\}=\{a,\{b,c,d\},e\}=\{a,b,\{c,d,e\}\}
\end{equation*} for all $a,b,c,d,e\in Z$. The importance of Abelian $JB^*$-triple systems derives from the fact that every $JB^*$-triple system is locally Abelian, which means that every element in a $JB^*$-triple system generates an Abelian subtriple. Every commutative $C^*$-algebra is an Abelian $JB^*$-triple system with the product $\{a,b,c\}=ab^*c$.
We call the elements of $$\Epi(Z,\C):=\{\vp:Z\to\C:\vp\neq0\text{ is a triple homomorphism}\}$$ the \textbf{characters}\index{characters} of Z. Following \cite{Kaup-ARiemannMappingTheoremForBoundedSymmetricDomainsInComplexBanachSpaces}, §1 we consider $\Epi(Z,\C)$ as a subspace of $Z'=B(Z,\C)$ and endow it with the $\sigma(Z^*,Z)$ topology. Then $\Epi(Z,\C)$ becomes a locally compact space and a principal $\mathbb{T}$-bundle for the group $\mathbb{T}=\{t\in\C:|t|=1\}.$
The base space $\Epi(Z,\C)/\mathbb{T}$ can be identified with the set of all $JB^*$-triple ideals $I\sub Z$ such that $Z/I$ is isometric to $\C$. The space
\[
    C_0^{\mathbb{T}}(\Epi(Z,\C)):=\{f\in C_0(\Epi(Z,\C))|\forall t\in\mathbb{T}\;\forall \lambda\in\Epi(Z,\C):f(t\lambda)=tf(\lambda) \}
\]
is a subtriple of the Abelian $C^*$-algebra $C_0(\Epi(Z,\C))$, the continuous functions on $\Epi(Z,\C)$ vanishing at infinity. The mapping
\begin{equation}\label{gelfand mapping}
    \hat{{}}:Z\to C_0^{\mathbb{T}}(\Epi(Z,\C))
\end{equation} defined by $\hat{x}(\lambda)=\lambda(x)$ for all $x\in Z$ and $\lambda\in\Epi(Z,\C)$ is called the \textbf{Gelfand transform}\index{Gelfand transform} of $Z$.

\begin{thm}[\cite{Kaup-OnJBtriplesdefinedbyfibrebundles}, Theorem 6.2]\label{abelsche trple equivalence}
For every $JB^*$-triple system $Z$ the following assertions are equivalent:
\begin{description}
\item[(a)] $Z$ is Abelian.
\item[(b)] $Z$ is a subtriple of a commutative $C^*$-algebra.
\item[(c)] The Gelfand transform of $Z$ is a surjective isometry onto $C_0^{\mathbb{T}}(\Epi(Z,\mathbb{C}))$.
\end{description}
\end{thm}
Especially, every Abelian $JB^*$-triple system is a TRO.

\begin{lemma}\label{abelsche triple sind TROs}
Let $Z$ be an Abelian $JC^*$-triple. Then $Z$ is a universally reversible TRO.
\begin{proof}
We only have to show that every Abelian $JC^*$-triple system is already a TRO since every TRO is already reversible, but by Theorem \ref{abelsche trple equivalence} we know that $Z$ is a subtriple of an Abelian $C^*$-algebra and therefore a TRO.
 \end{proof}
\end{lemma}

\begin{prop}\label{Abelian triple ternare hulle}
Let $Z$ be an Abelian $JC^*$-triple system, then \[
T^*(Z)\simeq C_0^{\mathbb{T}}(\Epi(Z,\mathbb{C}))
\] and the universal embedding $\rho_Z:Z\to C_0^{\mathbb{T}}(\Epi(Z,\mathbb{C}))$ is given by the Gelfand transform of $Z$.
\begin{proof} The Abelian $JC^*$-triple system $Z$ is by Lemma \ref{abelsche triple sind TROs} a universally reversible TRO.
Let $\hat{{}}:Z\to C_0^{\mathbb{T}}(\Epi(Z,\mathbb{C}))$ be the Gelfand transform, which is by Theorem \ref{abelsche trple equivalence} a $JB^*$-triple isomorphism. The identity mapping  $\id$ on $C_0^{\mathbb{T}}(\Epi(Z,\mathbb{C}))$ is, since we are in the Abelian world, also an antiautomorphism, satisfying $\id\circ\;\hat{{}}=\hat{{}}$. Since $\hat{Z}$ generates $C_0^{\mathbb{T}}(\Epi(Z,\mathbb{C}))$ as a TRO we obtain the statement from Lemma \ref{Tro identifizierer}.
\end{proof}
\end{prop}

\begin{defi}
Let $Z$ be an universally reversible $JC^*$-triple system. Define the \textbf{radical}\index{radical} of $Z$ to be the set $$\mathbf{R}(Z):=\bigcap_{\vp\in \Epi(Z,\C)\cup\{0\}}\ker(\vp).$$\end{defi}
In the case that $\Epi(Z,\C)=\emptyset$ we have $\ra(Z)=Z$.

The next proposition helps us to show that the radical of a universal reversible $JC^*$-triple system is universally reversible.

\begin{prop}
Let $Z$ be a universally reversible $JC^*$-triple system and $I\sub Z$ a $JB^*$-triple ideal, then $I$ is also universally reversible.
\begin{proof}
We assume that $T^*(I)\sub T^*(T)$.
It suffices to show that $\rho_Z(I)\sub T^*(Z)$ is reversible. Since $T^*(I)$ is a TRO-ideal and $\rho_Z(Z)$ is reversible by definition, we know that $\rho_Z(I)$ is reversible, if\[
\rho_Z(I)=T^*(I)\cap\rho_Z(Z).
\]
Let $x\in T^*(I)\cap\rho_Z(Z)$
 and $\pi:\rho_Z(Z)\to \rho_Z(Z)/\rho_Z(I)$ be the $JB^*$-quotient homomorphism. It follows from \cite{bohle-diss} Theorem 4.2.4 that $$T^*(Z)/T^*(I)\simeq T^*(\rho_Z(Z)/\rho_Z(I))$$ and therefore
$\pi(x)=\tau(\pi)(x)=0$, which yields $x\in\rho_Z(I)$.
\end{proof}
\end{prop}

Since the radical is always a $JB^*$-triple ideal the next corollary follows immediately.

\begin{cor}\label{radicalisunivrev}
  Let $Z$ be a universally reversible $JC^*$-triple system, then $\mathbf{R}(Z)$ is universally reversible.
 \end{cor}

\begin{thm}\label{gen of timon}
Let $T$ be a universally reversible TRO embedded in a $C^*$-algebra $\fA$ such that there exists a TRO antiautomorphism $\theta:\fA\to\fA$ of order $2$. Then we have an exact sequence of TROs
\begin{align}\label{radikalsequenz}
0\longrightarrow \ra(T)\oplus\theta(\ra(T))\longrightarrow  T^*(T)\longrightarrow C_0^{\mathbb{T}}\left(\Epi(T/\ra(T),\C)\right)\longrightarrow 0.
\end{align}
\begin{proof}
By Corollary \ref{radicalisunivrev} we know that the radical $\ra(T)$ is universally reversible and does not contain a TRO-ideal of codimension $1$ by construction. Using Lemma \ref{codimension Tro} we get $$ T^*(\ra(T))=\ra(T)\oplus\theta(\ra(T)).$$
The quotient $T/\ra(T)$ is an Abelian $JB^*$-triple system and we get with Proposition \ref{Abelian triple ternare hulle} that
$$T^*(T/(\ra(T)))=C_0^{\mathbb{T}}\left(\Epi(T/\ra(T),\C)\right).$$
The exactness of (\ref{radikalsequenz}) follows now from the exactness of
\begin{align*}
0\longrightarrow \ra(T)\longrightarrow  T\longrightarrow T/\ra(T)\longrightarrow 0,
\end{align*}
and \cite{bohle-diss} Theorem 4.2.4.
\end{proof}
\end{thm}

Theorem \ref{gen of timon} is a generalization of Lemma \ref{codimension Tro}. If we add the additional assumption that $T$ does not contain a one codimensional TRO-ideal, then $\ra(T)=T$ and thus \eqref{radikalsequenz} becomes
\begin{align*}
0\longrightarrow T\oplus\theta(T)\longrightarrow  T^*(T)\longrightarrow 0\longrightarrow 0.
\end{align*}

\bibliographystyle{alpha}

\bibliography{litarbeit}

\end{document}